\documentstyle[12pt]{article}
\begin{document}
\title{Beta and Gamma functions of Cayley-Dickson numbers}
\author{S.V. Ludkovsky}
\date{27.05.2004}
\maketitle

\section{Introduction.}
This paper continuous investigations of
function theory over Cayley-Dickson algebras \cite{luoyst,luoyst2}.
Cayley-Dickson algebras ${\cal A}_v$ over the field of real numbers
coincide with the field $\bf C$ of complex numbers for $v=1$,
with the skew field of quaternions, when $v=2$,
with the division nonassociative noncommutative algebra $\bf K$ of octonions
for $v=3$, for each $v\ge 4$ they are nonassociative and not
division algebras. The algebra ${\cal A}_{v+1}$ is obtained from
${\cal A}_v$ with the help of the doubling procedure.
This work provides examples of ${\cal A}_v$-meromorphic functions
and usages of line integrals over ${\cal A}_v$.
Here notations of previous papers \cite{luoyst,luoyst2} are used.
Discussions of references and results of others authors
can be found in \cite{luoyst,luoyst2} as well as physical applications
(see also \cite{baez,emch,guetze,hamilt,kansol,kurosh,lawmich,
rothe,ward} and references therein).
Beta and Gamma functions illustrate general theory of meromorphic
functions of Cayley-Dickson numbers and also applications of
line integration over ${\cal A}_v$.
\par  The results below show
some similarity with the complex case and as well differences
caused by noncommutativity and nonassociativity of Cayley-Dickson
algebras. It is necessary to mention that before works \cite{luoyst,luoyst2}
there was not any publication of others authors devoted to the line
integration of continuous functions of Cayley-Dickson numbers
or even quaternions along rectifiable paths. In works of others
authors integrations over submanifolds of codimension $1$
in $\bf H$ or $\bf K$ were used instead of line integral.
Therefore, in this respect publications \cite{luoyst,luoyst2} are
the first devoted to (integral) holomorphic functions of
Cayley-Dickson numbers.
\par If $g$ is a complex holomorphic function on a domain $V$ in the
complex plane $\Pi $ embedded into ${\cal A}_v$ and $g$
has a local expansion $g(z) = \sum_{n=0}^{\infty } a_n(z-z_0)^n$
in the ball $B(\Pi ,z_0,r^{-})$ for each $z_0\in Int (V)$,
where $Int (V)$ is the interior of
$V$ in $\Pi $, $r>0$ and $a_n=a_n(z_0)\in \bf C$ may depend on
parameter $z_0$,
$B(X,a,r^{-}):=\{ x\in X:$ $\rho (x,a)<r \} $ for a metrizable space
$X$ with metric $\rho $, $f$ is a function on a domain $U$ in
${\cal A}_v$, $v\ge 2$, such that $V\subset U\cap \Pi $ and
$f(z)=\sum_{n=0}^{\infty } a_n(z-z_0)^n$ is the local expansion of $f$
in $B({\cal A}_v,z_0,r^{-})$ for each $z_0\in V$, then
the line integral $\int_{\omega }f(z)dz$ over ${\cal A}_v$
for a rectifiable path $\omega $ in $V$ coincides with the classical
complex Cauchy line integral, since ${\hat f}|_V=f|_V$.
Nevertheless, if $f$ is an
${\cal A}_v$-holomorphic function on a domain $U$ in ${\cal A}_v$,
$V=U\cap \Pi $ is a domain in the complex plane $\Pi $ embedded
into ${\cal A}_v$, then in general $\int_{\omega }f(z)dz$
can not be reduced to Cauchy line integral for any rectifiable
path $\omega $ in $V$, since the generalized operator ${\hat f}$
is defined by values of $f$ in the neighbourhood of $\omega $
(see \cite{luoyst,luoyst2}).
\section{ Beta and Gamma functions of Cayley-Dickson numbers.}
\par {\bf 1. Definition.} The Gamma function is defined by the formula:
\par $(1)$ $\Gamma (z):=\int_0^{\infty }e^{-t}t^{z-1}dt$, \\
whenever this (Eulerian of the second kind)
integral converges and defined by ${\cal A}_v$-holomorphic
continuation elsewhere,
where $0<t\in \bf R$, $t^z:=exp (z ln t)$, $z\in {\cal A}_v$,
$t^z$ take its principal value, $dt$ corresponds to the Lebesgue
measure on $\bf R$, $ln : (0,\infty )\to \bf R$ is the classical
(natural) logarithmic function.
\par Denote by $Re (z):=(z+z^*)/2$ the real part of $z\in {\cal A}_v$,
${\cal I}_v:= \{ z\in {\cal A}_v: Re (z)=0 \} $, where $z^*$
is the conjugate of a Cayley-Dickson number $z$.
\par {\bf 2. Proposition.} {\it The gamma function has as singularities
only simple poles at the points $z\in \{ 0, -1, -2,... \} $
and $res (-n, \Gamma )M =[(-1)^n/n!]M$ for each $M\in {\cal I}_v$.}
\par {\bf Proof.} Write $\Gamma (z)$ in the form:
\par $(i)$ $\Gamma (z)=\Phi (z) + \Psi (z)$, where
\par $(ii)$ $\Phi (z)=\int_0^1e^{-t} t^{z-1}dt,$
\par $(iii)$ $\Psi (z):=\int_1^{\infty }e^{-t}t^{z-1}dt$. \\
Since $|e^{a+M}|=e^a$ for each $a\in \bf R$ and $M\in {\cal I}_v$
(see Corollary 3.3 \cite{luoyst2}), then $|t^{z-1}|\le t^{\delta -1}$
for each $Re (z)\le \delta $, where $\delta >0$ is a marked number.
From $\lim_{t\to \infty } e^{-t} t^{z-1}=0$ it follows, that
there exists $C=const >0$ such that $|e^{-t}t^{z-1}|\le Ce^{-t/2}$
for each $t>0$ and each $z$ with $Re (z)\le \delta $.
Therefore, $\Psi (z)$ is the ${\cal A}_v$-holomorphic function
in ${\cal A}_v$.
\par Consider change of variables $t=1/u$,
then $\Phi (z)=\int_1^{\infty }e^{-1/u}u^{-z-1}du$ for each $\delta >0$
and each $z$ with $Re (z)\le \delta $, hence
$|e^{-1/u} u^{-z-1}|\le u^{-\delta -1}$. Therefore, $\Phi (z)$ is
${\cal A}_v$-holomorphic, when $Re (z)>0$.
Substituting the Taylor series for $e^{-t}$ into the integral
expression $(ii)$, we get
\par $(iv)$ $\Phi (z)=\sum_{n=0}^{\infty }
(-1)^n \int_0^1t^{n+z-1}dt/n!=\sum_{n=0}^{\infty }(-1)^n(n+z)^{-1}/n!$.
Series $(iv)$ is uniformly and aboslutely convergent in any closed
domain in ${\cal A}_v\setminus \{ 0, -1, -2, ... \} $ and this series
gives ${\cal A}_v$-analytic continuation of $\Phi (z)$.
Thus $\Gamma (z)$ has only simple poles at the points
$z\in \{ 0, -1, -2, ... \} $.
\par The following Tannery lemma is true for ${\cal A}_v$-valued
functions (for complex valued functions see \S 9.2 \cite{copson}).
\par {\bf 3. Lemma.} {\it If $g(t)$ and $f(t,n)$ are functions
from $[a,\infty )$ to ${\cal A}_v$, $v\ge 2$, $\lim_{n\to \infty }
f(t,n)=g(t)$, $\lim_{n\to \infty }\lambda _n=\infty $, then
$\lim_{n\to \infty } \int_a^{\lambda _n}f(t,n)dt=\int_a^{\infty }
g(t)dt$, provided that $f(t,n)$ tends to $g(t)$ uniformly
on any fixed interval, and provided also that there exists a positive
function $M(t)$ such that $|f(t,n)|\le M(T)$ for each values of
$n$ and $t$ and such that $\int_a^{\infty }M(t)dt$ converges.}
\par {\bf Proof.} The sequence $f(t,n)$ converges uniformly to
$g(t)$ in the fixed segment $a\le t\le b$, $a<b<\infty $.
Using triangle inequalities and $| \int_a^bg(t)dt|\le \int_a^b
|g(t)|dt$ gives: $\limsup_{n\to \infty } | \int_a^{\lambda _n}
f(t,n)dt - \int_a^{\infty } g(t)dt| \le 2 \int_b^{\infty } M(t)dt$
for each $a<b<\infty $. From $\lim_{b\to \infty }
\int_b^{\infty } M(t)dt=0$ the statement of this lemma follows.
\par {\bf 4. Proposition.} {\it If $\Gamma (z,n) := n! n^z
[(...((z(z+1))(z+2))...(z+n)]^{-1}$, $n\in \bf N$, then
$\Gamma (z,n)$ tends to $\Gamma (z)$ as $n\to \infty $, the
convergence being uniform in any bounded canonical closed subset
$U\subset {\cal A}_v$ which contains no any of the singularities
of $\Gamma (z)$, $v\ge 2$.}
\par {\bf Proof.} Since ${\cal A}_v$ is power-associative and
$\bf R$ is the centre of the Cayley-Dickson algebra,
then $\{ n^z [z(z+1)(z+2)...(z+n)]^{-1} \} _{q(n+2)}$
does not depend on the order of multiplication regulated by
the vector $q(n+2)$ (see \cite{luoyst2}).
Therefore, $\Gamma (z,n)=(n/(n+1))^z z^{-1} \prod_{m=1}^n \{ (1+1/m)^z
(1+z/m)^{-1} \} $. Then $(1+1/m)^z(1+z/m)^{-1}=1+z(z-1)/(2m^2)+
O(1/m^3),$ when $m>0$ is large, hence $z^{-1} \prod_{m=1}^n \{ (1+1/m)^z
(1+z/m)^{-1} \} $ converges uniformly and aboslutely in any
bounded canonical closed domain $U$ in ${\cal A}_v$ to an
${\cal A}_v$-holomorphic function in accordance with
Theorem  3.21 \cite{luoyst2}. In view of Formulas $(3.6,7)$
in \cite{luoyst2} $|(1-t/n)^nt^{z-1}|=(1-t/n)^nt^{a-1}\le e^{-t}t^{a-1}$,
where $a:=Re (z)$. From $\int_0^n(1-t/n)^nt^{z-1}dt=n^z\int_0^1(1-u)^n
u^{z-1}du$ and integrating by parts we get
$\Gamma (z,n)=\int_0^n(1-t/n)^nt^{z-1}dt$, hence
$\lim_{n\to \infty }\Gamma (z,n)=\int_0^{\infty }e^{-t}
t^{z-1}dt=:\Gamma (z)$ by Lemma 3.
\par {\bf 5. Remark.} From the proof of Proposition 4 it follows,
that $\Gamma (z)=z^{-1}\prod_{m=1}^{\infty } \{ (1+1/m)^z
(1+z/m)^{-1} \} $ for each $z\in {\cal A}_v\setminus
\{ 0, -1, -2,... \} $, $v\ge 1$. The latter is known as the Euler's formula
in the case of complex numbers.
\par {\bf 6. Proposition.} {\it The gamma function satisfies identities:
\par $(i)$ $\Gamma (z+1)=z\Gamma (z)$ and
\par $(ii)$ $\Gamma (z) \Gamma (1-z)=\pi \csc (\pi z)$ \\
for each $z\in {\cal A}_v\setminus \{ 0, -1, -2,... \} $,
$v\ge 2$.}
\par {\bf Proof.} In view of power associativity of ${\cal A}_v$
and that $\bf R$ is the centre of the Cayley-Dickson algebra we get
\par $\Gamma (z+1)=z\lim_{n\to \infty } n!n^z[z(z+1)...(z+n)]^{-1}
n/(z+n+1)=z\Gamma (z)$, also
\par $\Gamma (z)\Gamma (1-z)=\lim_{n\to \infty }
\{ z(1-z^2/1^2)(1-z^2/2^2)...(1-z^2/n^2)
(1+(1-z)/n) \} ^{-1}=\{ z \prod_{n=1}^{\infty } (1-z^2/n^2) \} ^{-1}$.
In view of \S 6.83 $\pi \csc (\pi z)=
\{ z \prod_{n=1}^{\infty } (1-z^2/n^2) \} ^{-1}$ for complex
$z$ in ${\bf C}\setminus \bf Z$ \cite{copson}.
Using the ${\cal A}_v$-holomorphic
extension of this function from the complex domain onto the
corresponding domain ${\cal A}_v\setminus \bf Z$
(see Proposition 3.13, Corollary 2.13 and Theorems 3.10,
3.21 \cite{luoyst2}), we get Formula $(ii)$.
\par {\bf 7. Definition.} A function $F$ on an unbounded domain
$U$ in ${\cal A}_v$, $v\ge 2$, is said to have an asymptotic
expansion $F\sim \sum_{|k|\le 0}(a_k,z^k)$, if
\par $\lim_{z\in U, |z|\to \infty }
z^n \{ F(z) - \sum_{|k|\le 0}(a_k,z^k) \} =0$ \\
for each $n\in \bf N$, where
$k = (k_1,...,k_n),$ $|k| := k_1+...+k_n$, $k_j\in \bf Z$
for each $j$, $n\in \bf N$, $(a_k,z^k):=a_{k_1}z^{k_1}...a_{k_n}z^{k_n}$,
$a_{k_j}\in {\cal A}_v$ for each $j$.
\par We write $F(z)\sim G(z) \sum_{|k|\le 0}(a_k,z^k)$, if
$G(z)^{-1} F(z) \sim \sum_{|k|\le 0}(a_k,z^k)$.
The term $G(z) a_0$ is called the dominant term of the asymptotic
representation of $F(z)$.
\par {\bf 8. Lemma.} {\it Let $f(t)$ be a function
in an unbounded domain $U$ in ${\cal A}_v$ possibly with a
branch point at $0$ and such that \\
$f(z) = \sum_{m=1}^{\infty } a_mz^{(m/r)-1} $, \\
when $|z|\le a$, $a>0$, $r>0$, let also $f$ be ${\cal A}_v$-holomorphic
in $B(U,0,a+\delta )\setminus \{ 0 \} $, where $\delta >0$.
Suppose, that when $t\ge 0$, $|f(t)|< C e^{bt}$, where
$C>0$ and $b>0$ are constants. Then
\par $F(z)=\int_0^{\infty } e^{-zt}f(t)dt \sim
\sum_{n=1}^{\infty } a_n \Gamma (n/r) z^{-n/r}$, \\
when $|z|$ is large and $|Arg (z)|\le \pi /2 - \epsilon $,
where $\epsilon >0$ is arbitrary.}
\par {\bf Proof.} For each $n\in \bf N$ there exists a constant
$C=const >0$ such that
\par $|f(t) - \sum_{m=1}^{n-1} a_mt^{(m/r)-1}|\le C t^{(n/r)-1}e^{bt}$ \\
for each $t\ge 0$. In view of Formulas
$(3.2,3)$ \cite{luoyst2}
\par $|\int_0^{\infty }e^{-zt} [f(t)-\sum_{m=1}^{n-1}
a_mt^{(m/r)-1}]dt|\le \int_0^{\infty }e^{-xt}Ct^{(n/r) -1}e^{bt}dt$ \\
$= C\Gamma (n/r) (x-b)^{-n/r} $ \\
for each $x>b$, where $x:=Re (z)$.
From the condition $| Arg (z)|\le \pi /2 -\epsilon $ it follows,
that $x\ge |z| \sin (\epsilon )$, such that $x>b$ for
$|z|>b\csc (\epsilon )$.  Therefore, for $| Arg (z)|\le \pi /2-\epsilon
<\pi /2$ and $|z|>b \csc (\epsilon ) ,$ there is the inequality: \\
$|z^{n/r} \int_0^{\infty }e^{-zt}[f(t)-\sum_{m=1}^{n-1}a_mt^{(m/r)-1}]
dt|\le C \Gamma (n/r) |z|^{n/r}/(|z| \sin (\epsilon ) - b)^{n/r} = O(1)$.
\par {\bf 9. Proposition.} {\it Let $0<\delta <\pi /2$,
$z\in {\cal A}_v\setminus \{ 0, -1, -2,...\} $,
$|Arg (z)| \le \pi - \delta $, $v\ge 2$. Then there exists
the asymptotic expansion:
\par $Ln \Gamma (z) \sim (z-1/2) Ln (z) - z + (ln (2\pi ))/2
+ \sum_{n=1}^{\infty } (-1)^{n-1}B_n [2n(2n-1)z^{2n-1}]^{-1}$, \\
where $B_n$ are Bernoulli numbers defined by the equation:
$(z/2) \coth (z/2) = 1+\sum_{n=1}^{\infty }(-1)^{n-1} B_nz^{2n}/(2n)!$.}
\par {\bf Proof.} If $z> 0$, then the substitution $t=zu$
gives $\Gamma (z)=\Gamma (1+z)/z = z^ze^{-z}\int_0^{\infty }
(ue^{1-u})^zdu$ and by analytic continuation the formula
$\Gamma (z)=z^ze^{-z}\int_0^{\infty } (ue^{1-u})^zdu$
is true for each copy of $\bf C$, $0\in \bf C$, embedded into ${\cal A}_v$.
In view of independence of this formula from such embedding
and power associativity of ${\cal A}_v$ it follows, that
it is true for each $z\in {\cal A}_v$ with $Re (z) > 0$.
For $Re (z)>0$ and large $|z|$ using substitutions
$e^{-t}=\eta e^{1-\eta }$ for $t\in (0,\infty )$ and
$\eta \in (1,\infty )$; also the substitution
$e^{-t}=ue^{1-u}$ for $t$ decreasing monotonously from $\infty $
to $0$ and $u\in (0,1)$, we get
$z^{-z}e^z\Gamma (z)=\int_0^{\infty } e^{-zt} (d\eta /dt -
du/dt)dt$. Consider two real solutions $\eta $ and $u$ of the equation
$t=u-1-ln (u)$ and the equation $\zeta ^2 /2 =w-Ln (1+w)$
which defines $w=w(\zeta )$ for $\zeta \in {\bf R}\oplus M\bf R$.
It has two branches $\zeta =\beta w (1-2w/3+2w^2/4-...)^{1/2}$,
where $\beta = -1$ or $\beta =1$. Each branch is the analytic function
of $w$ in the domain $ \{ w\in {\bf R}\oplus M{\bf R}:$
$|w|<1 \} $ with a simple zero at $w=0$. For $\beta =1$ there exists
a unique solution $w=\zeta +a_2\zeta ^2+a_3\zeta ^3+...$
in $ \{ \zeta : |\zeta |<\rho \} $, $na_nM=res (0, \zeta ^{-n})M$
for each $n>1$.
\par Thus $w$ has two branches $w_1$ and $w_2 (\zeta )=
w_1 (-\zeta )$. Singularities of $w(\zeta )$ are only points
at which $dw/d\zeta $ is zero or infinite, hence these are
$\zeta =0$, also points corresponding to $w=0$ and $w=-1$, since
$dw/d\zeta =\zeta (1+w)/w$. Then $\zeta =0$ is not a branch-point
of $w_1$, to $w=-1$ there corresponds $\zeta =\infty $.
Therefore, singularities are: $\zeta ^2=4n\pi M$, where
$n\in {\bf Z}\setminus \{ 0 \} $. Then $\eta $ and $u$ are
${\cal A}_v$-holomorphic, when $|(z+{\tilde z})/2|<2\pi $
possibly besides $z=0$ and when $|z|<2\pi $, where
$\zeta ^2 =: 2z$ such that $\eta = 1+(2z)^{1/2}+a_2(2z)+a_3(2z)^{3/2}
+a_4(2z)^2+...$, $u = 1-(2z)^{1/2}+a_2(2z)-a_3(2z)^{3/2}+a_4(2z)^2- ...$,
the square roots are taken positive, when $z>0$.
Applying Lemma 8 we get the asympotic expansion.
In view of Theorem 2.15 \cite{luoyst2}
for $M\in {\cal I}_v$ with $|M|=1$ and $\alpha \in \bf R$
and a loop defined by $\rho e^{Mt}$ on the boundary of the
sector $|z|\le \rho $ and two lines $Arg (z)=0$,
$Arg (z)=M\alpha $, where $\alpha \in (- \pi /2, \pi /2)$,
$g(z):=d(\eta -u)/dt$, provides the equality: $\int_0^{\infty } e^{-zt}
g(t)dt=\int_0^{\infty } \exp ( - zte^{M\alpha } ) g(te^{M\alpha })
e^{M\alpha } dt$, when $Arg (z)\in M\bf R$,
$Re (z)>0$ and $Re (ze^{M\alpha })>0$,
since ${\bf R}\oplus M\bf R$ is isomorphic with $\bf C$
which is commutative. Therefore, the latter integral converges
uniformly and provides the analytic function.
Two regions $Re (z)>0$ and $Re (ze^{M\alpha })>0$ have
a common area and by the analytic continuation:
$z^{-z}e^z \Gamma (z)=\int_0^{\infty } \exp (-zt e^{M\alpha })
g(te^{M\alpha }) e^{M\alpha }dt$, when
$\alpha \in (-\pi /2, \pi /2)$. Applying Lemma 8 we get the region
of validity of this asymptotic expansion, since $M$ is arbitrary.
\par {\bf 10. Corollary.} {\it For large $|y|$ there is the asymptotic
expansion $|\Gamma (x+My)|\sim (2\pi )^{1/2} |y|^{x-1/2}\exp
(-\pi |y|/2)$ uniformly by $M\in {\cal I}_v$, $|M|=1$, where
$v\ge 2$, $y\in \bf R$.}
\par {\bf 11. Corollary.} {\it $\pi ^{1/2} \Gamma (2z)=2^{2z-1}
\Gamma (z) \Gamma (z+1/2)$ for each $z\in {\cal A}_v
\setminus \{ 0, -1, -2, ... \} $, $v\ge 2$.}
\par The proof is analogous to \S \S 9.55, 9.56 \cite{copson},
since ${\bf R}\oplus M\bf R$ is isomorphic with $\bf C$
for each $M\in {\cal I}_v$, $v\ge 2$, $|M|=1$.
\par {\bf 12. Proposition.} {\it For all $z\in {\cal A}_v$:
\par $1/\Gamma (z)=(2\pi )^{-1}
(\int_{\psi } e^{\zeta }\zeta ^{-z}d\zeta )M^*$ \\
for a loop $\psi $ and $z$ in the plane ${\bf R}\oplus M\bf R$,
$M\in {\cal I}_v$, $|M|=1$, $\psi $ starts at $ - \infty $ of
the real axis, encircles $0$ once in the positive direction and returns
to the starting point.}
\par {\bf Proof.} Consider the integral
$\int_{\psi }e^{\zeta }\zeta ^{-z}d\zeta =: \int_{\psi } f(\zeta )d\zeta $,
the integrand $f(\zeta )$ 
has a branch point at zero, but each branch is a one-valued function
of $\zeta $ and each branch is ${\cal A}_v$-holomorphic in
${\cal A}_v\setminus Q$, where $Q$ is a submanifold in ${\cal A}_v$
of real codimension $1$ such that $(-\infty ,0]\subset Q$
(see \S 3.7 \cite{luoyst}). Then take a branch $e^{\zeta }\zeta ^{-z}
=\exp (\zeta -z Ln (\zeta ))$, where $Ln (\zeta )$ takes its principal
value. Consider a rectifiable loop $\gamma $ in ${\bf R}\oplus
M\bf R$ encompassing zero in the positive direction and beginning
at $- \rho $ on the lower edge of the cut and returns to $- \rho $
at the upper edge of the cut, where $\rho >0$.
\par In view of Theorem
2.15  \cite{luoyst2} the value of the integral
is not changed by the deformation to a contour $\gamma $ consisting
of the lower edge of the cut intersected with $[- \rho , - \delta ]$,
where $0<\delta <\rho $, the circle $|z|=\delta $ in the plane
${\bf R}\oplus M\bf R$, and the upper edge of the cut intersected with
$[- \rho , - \delta ]$. On the upper edge of the cut in $\gamma $:
$\zeta =ue^{\pi M}$, where $u>0$, $u\in \bf R$, and
$f(\zeta )=\exp (-u-z ln (u) -z\pi M)=e^{-u}u^{-z}e^{-z \pi M}$.
On the lower egdge of the cut in $\gamma $:
$\zeta =ue^{-\pi M}$ and $f(\zeta )=e^{-u}u^{-z}e^{z\pi M}$.
Therefore, $\int_{\gamma }f(\zeta )d\zeta = (e^{z\pi M}-e^{-z\pi M})
\int_{\delta }^{\rho }e^{-u}u^{-z}du+J$, where $J:=
\int_{-\pi }^{\pi }\exp (\delta e^{\theta M})M
\delta ^{(1-z)} e^{(1-z) \theta M}d\theta $, since ${\bf R}\oplus M\bf R$
is isomorphic with $\bf C$ and ${\hat f}(z)|_{{\bf R}\oplus M\bf R}.h=
f(z)h$ for each $h$ and $z\in {\bf R}\oplus M\bf R$ (see
Theorem 2.7 \cite{luoyst2}).
\par If $z=x+yM$, where $x$ and $y\in \bf R$,
then $|J|\le \int_{- \pi }^{\pi } \delta ^{1-x} \exp (\delta \cos (\theta )
+y\theta )d\theta \le 2\pi \delta ^{1-x} e^{\delta +\pi |y|}$,
consequently, $\lim_{\delta \to 0} J=0$, when $x<1$.
Hence $\int_{\gamma }e^{-\zeta }\zeta ^{-z}d\zeta =
2\sin (\pi z) (\int_0^{\rho }e^{-u}u^{-z}dz) M$ for $Re (1-z)>0$.
Suppose $\psi $ is the loop obtained from $\gamma $ by tending
$\rho $ to the infinity, then $\int_{\psi }e^{\zeta }\zeta ^{-z}d\zeta
=2\sin (\pi z) (\int_0^{\infty }e^{-u}u^{-z}du)M=
2\sin (\pi z) \Gamma (1-z)M$,
since $\Gamma (z)\Gamma (1-z)=\pi \csc (\pi z)$, hence
$1/ \Gamma (z)= (2\pi )^{-1}(\int_{\psi }e^{\zeta }\zeta ^{-z}d\zeta )M^*$.
Since $M\in {\cal I}_v$ with $|M|=1$ is arbitrary, then this formula
is true in ${\cal A}_v\setminus Q$ with $Re (1-z)>0$. By the
complex holomorphic continuation this formula is true for all
values of $z$ in ${\bf R}\oplus M\bf R$.
\par {\bf 13. Corollary.} {\it Let $M$ and $\psi $ be as in
Proposition 12, then $\Gamma (z)=(2\sin (\pi z))^{-1}
(\int_{\psi }e^{\zeta }\zeta ^{z-1}d\zeta )M^*$
for each $z\in {\bf R}\oplus M{\bf R}\setminus \bf Z$.}
\par {\bf 14. Definition.} The Beta function $B(p,q)$ of Cayley-Dickson
numbers $p, q\in {\cal A}_v$, $v\ge 2$, is defined by the equation:
\par $B(p,q):=\int_0^1\zeta ^{p-1}(1-\zeta )^{q-1}d\zeta $, \\
whenever this integral (Eulerian of the first kind) converges,
where $\zeta ^{p-1}:=e^{(p-1) ln (\zeta ) }$ and the logarithm
has its principal value. This equation defines $B(p,q)$
for each $Re (p)>0$ and $Re (q)>0$. For others values of
$p$ and $q$ it is defined by the complex holomorphic continuation
by $p$ and $q$ separately and subsequently in each complex plane
${\bf R}\oplus M\bf R$ and ${\bf R}\oplus S\bf R$, $M, S\in {\cal I}_v$,
$|M|=1$ and $|S|=1$. 
\par {\bf 15. Proposition.} {\it Let $p, q\in {\cal A}_v$,
$v\ge 2$, such that the minimal subalgebra $\Upsilon _{p,q}$
containing $p$ and $q$ has embedding into $\bf K$, then \\
$B(p,q) - B(q,p) = [ B(p,q) - B(p,q_0-q') - B(p_0-p',q) +
B(p_0-p',q_0-q')](q')^*{q'}_2/2$,  \\
where $p_0 := Re (p)$, $p':= p - Re (p)$, ${q'}_2\perp p'$,
${q'}_1\parallel p'$ relative to the scalar product
$(z,\eta ):=Re(z\eta ^*)$, $q'={q'}_1+{q'}_2$.}
\par {\bf Proof.} Making the substitution
$\eta \mapsto 1 - \eta $ of the variable, we get
$\int_0^1\eta ^{q-1}(1-\eta )^{p-1}d\eta =\int_0^1(1-\eta )^{q-1}
\eta ^{p-1}d\eta $, but in general $p$ and $q$ do not commute.
In view of Formulas $(3.2, 3.3)$ \cite{luoyst2} the commutator
of two terms in the integral is: \\
$[t^{p-1}, (1-t)^{q-1}] = 2 t^{p_0-1} (1-t)^{q_0-1}
[(\sin |p'ln t|)/|p' ln t|) (\sin |q' ln (1-t)|)/|q' ln (1-t)|)]
(p' ln t) (q_2' ln (1-t)).$ \\
On the other hand, $[(\sin |M|) / |M|]M = [e^M-e^{-M}]/2$
for each $M\in {\cal I}_v$, hence  \\
$\int_0^1[t^{p-1},(1-t)^{q-1}]dt= (\int_0^1t^{p_0-1}(1-t)^{q_0-1}
[t^{p'}-t^{-p'}] [(1-t)^{q'}-(1-t)^{-q'}]dt) (q')^*{q'}_2/2$ \\
$= [ B(p,q) - B(p,q_0-q') - B(p_0-p',q) +
B(p_0-p',q_0-q')](q')^*{q'}_2/2$, since $\bf K$ is alternative and
$p'{q'}_2=p'((q'{q'}^*){q'}_2)=p'(q'({q'}^*{q'}_2))$.
\par {\bf 16. Remark.} Let $G$ be a classical Lie group over $\bf R$
and $g = T_eG$ be its Lie algebra (finite dimensional over $\bf R$).
Suppose that $e: V\to U$ is the exponential mapping of the neighbourhood
$V$ of zero in $g$ into a neighbourhood $U$ of the unit element
$e\in G$, $ln: U\to V$ is the logarithmic mapping. 
Then $w = ln (e^u\circ e^v)$, $w = w(u,v)$, is given
by the Campbell-Hausdorff formula
in terms of the adjoint representation $(ad\quad u)(v):=[u,v]$:
$$w=\sum_{n=1}^{\infty }n^{-1}\sum_{r+s=n, r\ge 0, s\ge 0}({w'}_{r,s}
+{w"}_{r,s}),\mbox{ where } {w'}_{r,s}=\sum_{m=1}^{\infty }$$
$$(-1)^{m-1}m^{-1}\sum^*((\prod_{i=1}^{m-1}(ad\quad u)^{r_i}(ad\quad v)^{s_i}
(r_i!)^{-1}(s_i!)^{-1})(ad \quad u)^{r_m}(r_m!)^{-1})(v),$$
$${w"}_{r,s}=\sum_{m=1}^{\infty }
(-1)^{m-1}m^{-1}\sum^{**}((\prod_{i=1}^{m-1}(ad\quad u)^{r_i}
(ad\quad v)^{s_i}(r_i!)^{-1}(s_i!)^{-1})(u),$$
$\sum^*$ means the sum by $r_1+...+r_m=r,$ 
$s_1+...+s_{m-1}=s-1,$ $r_1+s_1\ge 1,$ ..., $r_{m-1}+s_{m-1}\ge 1,$ 
$\sum^{**}$ means the sum by $r_1+...+r_{m-1}=r-1,$
$s_1+...+s_{m-1}=s,$ $r_1+s_1\ge 1,$ ..., $r_{m-1}+s_{m-1}\ge 1.$
In particular, this formula can be applied to the multiplicative
group $G={\bf H}\setminus \{ 0 \} $ with $U=G$ and $V=g$, since each
quaternion can be represented as a $2\times 2$ complex matrix,
where generators of $\bf H$ are Pauli matrices \cite{boug}.
\par {\bf 17. Theorem.} {\it Let $p, q \in {\cal A}_v$, $v\ge 2$,
such that the minimal subalgebra $\Upsilon _{p,q}$ generated
by $p$ and $q$ has embedding into $\bf H$, then \\
$\Gamma (p) \Gamma (q) = \Gamma (w(p,q)) B(p,q) - $  \\
$[\Gamma (w(p,q)) - \Gamma (w(p,q_0-q'))]{q'}^*{q'}_2
[B(p,q)-B(p_0-p',q)]/2$, \\
where $p_0 := Re (p)$, $p':= p - Re (p)$, ${q'}_2\perp p'$,
${q'}_1\parallel p'$ relative to the scalar product
$(z,\eta ):=Re(z\eta ^*)$, $q'={q'}_1+{q'}_2$, $w(p,q)$
is given in Remark 16.}
\par {\bf Proof.} Let $S_R:= \{ (x,y) \in {\bf R^2}: $
$0\le x\le R, 0\le y\le R \} $. Then
\par $\Gamma (p) \Gamma (q) =
\int_0^{\infty }e^{-x}x^{p-1}dx \int_0^{\infty } e^{-y}y^{q-1}dy$
\par $=\lim_{R\to \infty }\int_0^R( \int_0^R e^{-x-y}x^{p-1}
y^{q-1} dy)dx$
\par $=\lim_{R\to \infty }\int \int_{S_R}e^{-x-y}x^{p-1}y^{q-1}dx dy$ \\
in accordance with the Fubini theorem, since each function
$f: U\to \bf H$ has the form $f(z)=f_1(z)+f_i(z)i+f_j(z)j+f_k(z)k$
for each $z$ in a domain $U$ in $\bf H^n$, $f_1$, $f_i$, $f_j$, $f_k$
are real-valued functions, $ \{ 1,i,j,k \} $ are generators of $\bf H$.
Consider a triangle
$T_R := \{ (x,y)\in {\bf R^2}:$ $0\le x, 0\le y, x+y\le R \} $
and put $f(x,y) := e^{-x-y}x^{p-1}y^{q-1}$, where $p, q\in {\cal A}_v$ are
marked, then \\
$| \int \int_{S_R}f(x,y)dxdy - \int \int_{T_R}f(x,y)dxdy|
\le \int \int_{S_R\setminus T_R}|f(x,y)|dxdy$
\par $\le \int \int_{S_R}|f(x,y)|dxdy - \int \int_{S_{R/2}}|f(x,y)|dxdy$. \\
We have $\lim_{R\to \infty } \int \int_{S_R}|f(x,y)|dxdy =
\Gamma (p_0)\Gamma (q_0)$, hence \\
$\lim_{R\to \infty }\int \int_{S_R\setminus S_{R/2}}|f(x,y)|dxdy=0$.
Therefore, \\
$\Gamma (p)\Gamma (q) = \lim_{R\to \infty } \int \int_{T_R}e^{-x-y}
x^{p-1}y^{q-1}dxdy$. \\
The substitution $x+y=\xi $, $y=\xi \eta $ and application of the Fubini
theorem gives \\
$\Gamma (p)\Gamma (q)=\int_0^{\infty } \int_0^1 e^{-\xi }\xi ^p
(1-\eta )^{p-1}\xi ^{q-1}\eta ^{q-1}d\xi d\eta $, \\
since $\bf H$ is associative, $\xi ^{p-1}$ commutes with $(1-\eta )^{p-1}$,
$\xi ^{q-1}$ commutes with $\eta ^{q-1}$. Therefore, \\
$\Gamma (p)\Gamma (q)=\Gamma (w(p,q)) B(p,q)+
\int_0^{\infty }\int_0^1e^{-\xi }\xi ^p[(1-\eta )^{p-1},\xi ^{q-1}]
\eta ^{q-1}d\xi d\eta $. \\
Let $M$ and $N$ be in ${\cal I}_v$, then $e^Ne^M=(\cos |M|) e^N+
[(\sin |N|)/|N|]Me^{N_1-N_2}$, where $M\perp N_2$, $M\parallel N_1$
relative to the scalar product $(z, \eta ) := Re (z \eta ^*)$,
$N_1, N_2\in {\cal I}_v$, $N=N_1+N_2$ (see Formulas
$(3.2, 3.3)$ \cite{luoyst2}).
Therefore, \\
$\int_0^{\infty }\int_0^1e^{-\xi }\xi ^p[(1-\eta )^{p-1},\xi ^{q-1}]
\eta ^{q-1}d\xi d\eta =$ \\
$ - \int_0^{\infty } \int_0^1 e^{-\xi }\xi ^{p-1}
[\xi ^q - \xi ^{q_0-q'}]({q'}^*{q'}_2)[(1-\eta )^{p-1} -
(1-\eta )^{p_0-p'-1}] \eta ^{q-1}d\xi d\eta /2 $ \\
$=-[\Gamma (w(p,q)) - \Gamma (w(p,q_0-q'))]({q'}^*{q'}_2)
[B(p,q)-B(p_0-p',q)]/2$.
\par {\bf 18. Note.} Proposition 15 and Theorem 17 show differences
in identities for Beta and Gamma functions between commutative case of
$\bf C$ and noncommutative cases of ${\cal A}_v$, $v\ge 2$, and $\bf H$
particularly. Certainly, in the particular case if $\Upsilon _{p,q}$
has embedding into $\bf C$, then ${q'}^*{q'}_2=0$ and
Proposition 15 and Theorem 17 give classical results,
but for general $p$ and $q$ the subalgebra $\Upsilon _{p,q}$
can have no any embedding into $\bf C$.

\thanks{
Address: Sergey V. Ludkovsky, Mathematical Department, TW-WISK,
Brussels University, V.U.B., Pleinlaan 2, Brussels 1050, Belgium. \\
{\underline {Acknowledgment}}. The author thanks the Flemish
Science Foundation for support through the Noncommutative Geometry
from Algebra to Physics project.}
\end{document}